%% file: OvertwistedOpenBooksAndStallingsTwist.tex
\documentclass[]{%
amsart}

\input{OTOB-preamble}

\begin{document}


\title[overtwisted open books and Stallings twist]%
{open books supporting overtwisted contact
structures and Stallings twist}
\author{Ryosuke YAMAMOTO}
\date{}


\begin{abstract}
 We study open books (or open book decompositions)
 of a closed oriented $3$-manifold
 which support overtwisted contact
 structures.
 We focus on a simple closed curve along which one can perform
 Stallings twist, called ``twisting loop''.
 We show that
 the existence of a twisting loop on the fiber surface of an open book
 is equivalent up to positive stabilization to 
 the existence of an overtwisted disk in the contact manifold
 supported by the open book.
 We also show a criterion for overtwistedness using a certain arc properly
 embedded in the fiber surface, which is an extension of Goodman's one.
\end{abstract}

\maketitle

 \section{Introduction}\label{sec:intro}
 \input{intro}

 \section*{Acknowledgement}
 This paper is supported by 21st century COE program
 ``Constitution of wide-angle mathematical basis focused on knots''.

 \section{Preliminaries}\label{sec:prelim}
 \input{prelim}

 \section{Overtwisted open books}\label{sec:otob}
 \input{otob}


\bibliographystyle{plain}
%


\bigskip\noindent
Osaka City University Advanced Mathematical Institute \\
3-3-138, Sugimoto, Sumiyoshi-ku, Osaka, 558-8585 Japan

\bigskip\noindent
Email: ryosuke@sci.osaka-cu.ac.jp

\end{document}

%% file: OTOB-preamble.tex
\usepackage{amsmath}
\usepackage{amssymb}
\usepackage{amsthm}
\usepackage{amscd}

\usepackage[dvips]{graphicx}

\usepackage{extarrows}

\newcommand{\ud}{\mathrm{d}}

\newtheorem{thm}{Theorem}[section]
\newtheorem{lem}[thm]{Lemma}

\newtheorem{prop}[thm]{Proposition}

\theoremstyle{definition}
\newtheorem{dfn}{Definition}[section]
\newtheorem{rmk}[dfn]{Remark}



%% file: intro.tex

Stallings \cite{St} introduced two operations
which create a new open book (or open book decomposition)
of a closed oriented $3$-manifold
from another of the manifold.
One of them is called {\itshape Stallings twist\/},
which is a Dehn twist along a certain simple closed curve,
called a {\itshape twisting loop\/},
on the fiber surface of an open book.
The other is positive (negative resp.) stabilization
of an open book,
which is also known as plumbing of positive (negative resp.)
Hopf band to the fiber surface of an open book
(See Section \ref{sec:prelim}).
On $S^3$, Harer \cite{Ha} showed that every open book can be obtained
from the standard open book of $S^3$,
i.e., the open book with a $2$-disk as the fiber
and the identity map as the monodromy map,
by (de-)stabilizations and Stallings twists.
Moreover he conjectured that Stallings twists can be omitted.
This conjecture has been proved affirmatively
by Giroux \cite{Gi}.
He showed a one-to-one correspondence between
isotopy classes of contact structures on $M$ and
equivalence classes of open books on $M$
modulo positive stabilization (See \cite{GG} for further information).
We may say about this result that
a contact structure leads to a topological property of open books
via Giroux's one-to-one correspondence.
In this paper we will deal with a study in the opposite direction.
We show that
a twisting loop on the fiber surface of an open book is
related directly to an overtwisted disk in the contact structure
which is corresponding to the open book via Giroux's one-to-one
correspondence.

Let $M$ be a closed oriented $3$-manifold.
%
Denote by $(\Sigma,\phi)$ an open-book of $M$,
where $\Sigma$ is a fiber surface embedded in $M$
and $\phi$ is a monodromy map.
%
We say that an contact structure $\xi$ on $M$
is {\itshape supported}
by an open book $(\Sigma,\phi)$
if $(\Sigma,\phi)$ is a corresponding one to $\xi$
on Giroux's one to one correspondence.
We will investigate open books
supporting an overtwisted contact structure.
For simplicity, we call such open books
{\itshape overtwisted open books\/}.
We show the following:

\begin{thm}
 \label{thm:main}
 Let $(\Sigma,\phi)$ be an open book of a closed
 oriented $3$-manifold.
 The following are equivalent;
 \begin{enumerate}
  \item $(\Sigma,\phi)$ is overtwisted.
  \item $(\Sigma,\phi)$ is equivalent up to positive stabilization
	to an open book
	whose fiber surface has a twisting loop. 
  \item $(\Sigma,\phi)$ is equivalent up to positive stabilization
	to an open book	$(\Sigma',\phi')$ with an arc $a$
	properly embedded in $\Sigma'$ such that
	$i_{\partial}(a,\phi'(a)) \leq 0$. 
 \end{enumerate}
\end{thm}

The arc $a$ in $(3)$ is an extension of Goodman's {\itshape sobering
arc\/} \cite{Go}.
The boundary intersection number $i_{\partial}$ of $a$ and $\phi(a)$,
introduced by Goodman,
is defined in Section \ref{sec:otob}.
We should mention that
the equivalence between $(1)$ and $(3)$ has already shown by 
Honda, Kazez and Mati\'{c} \cite[Theorem 1.1]{HKM},
but the proof given here differs from theirs
in focusing a twisting loop to detect an overtwisted disk.
%


%% file: prelim.tex

 Let $M$ be a closed oriented $3$-manifold.
 We denote by $E(X)$ the exterior of $X$ in $M$
 and by $N(X)$ a regular neighbourhood of $X$ in $M$,
 where $X$ is a submanifold in $M$.

 \subsection{Open book}
 Let $K$ be a fibered knot or link in $M$,
 i.e., there is a fibration map
 $f:E(K)\rightarrow S^{1}={\mathbb R}/{\mathbb Z}$
 such that $f$ maps meridian of $K$ to $S^{1}$ homeomorphically.
 We denote by $\Sigma_{t}$ the fiber surfaces $f^{-1}(t)$
 for each $t\in {\mathbb R}/{\mathbb Z}$,
 and by $\Sigma$ the homeomorphism type of the fiber surface.
 We often identify the abstract $\Sigma$ and $\Sigma_{0}$ embedded in
 $M$.
 In this situation, $M$ has a decomposition as follows;
 \[
  M
  = (\Sigma\times [0,1]/(x,1)\sim(\phi(x),0))
  \cup_{g}
  (D^{2}\times\partial\Sigma),
 \]
 where $\phi$ is an automorphism of $\Sigma$
 fixing $\partial\Sigma$ pointwise,
 and $g$ is a gluing map between the boundary tori such that
 $g(\{p\}\times [0,1]/(p,1)\sim(p,0))
 =\partial(D^{2}\times\{p\})$
 for $p\in\partial\Sigma$.
 We call this structure of $M$ an {\itshape open book\/} of $M$
 and denote by a pair $(\Sigma,\phi)$.
 The automorphism $\phi$ is called a {\itshape monodromy map\/}
 of the open book.

 Let $c$ be a simple closed curve on $\Sigma$. 
 We use notation $\operatorname{Fr}(c; \Sigma)$
 for the framing of $c$ determined by a curve parallel to $c$ on $\Sigma$,
 and $D(c)$ ($D(c)^{-1}$ resp.)
 for a positive (negative resp.) Dehn twist on $\Sigma$ along $c$.
 We say that $c$ is essential on $\Sigma$
 if $c$ does not bound a disk region on $\Sigma$.

\begin{dfn}
 An essential simple closed curve $c$ on $\Sigma$
 is a {\itshape twisting loop}
 if $c$ bounds a disk $D$ embedded in $M$ 
 and satisfies that
 $\operatorname{Fr}(c; \Sigma)=\operatorname{Fr}(c; D)$.
\end{dfn}

\begin{figure}[ht]
 {\unitlength=1cm
 \begin{center}
  \begin{picture}(8.6,3.5)(0,0)
   \put(0,.6){\includegraphics[height=3cm]{%
   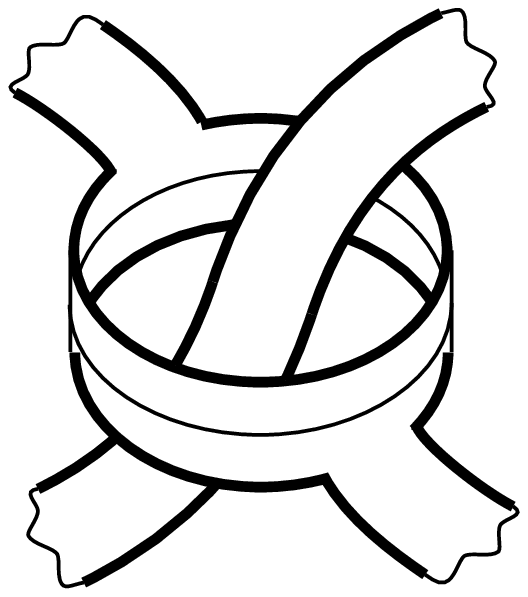} }       
   \put(5.6,.6){\includegraphics[height=3cm]{%
   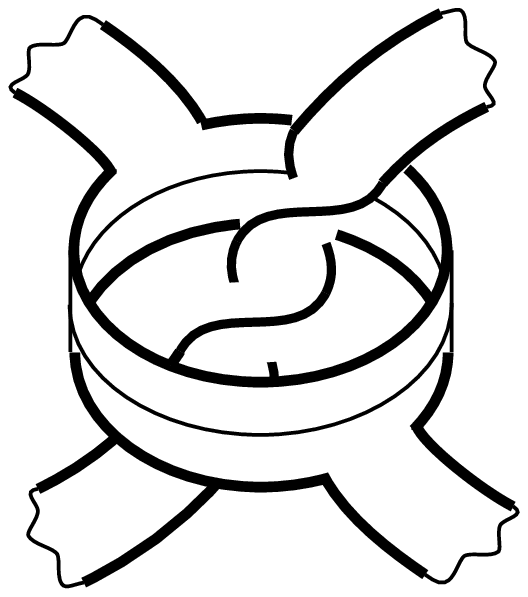} }
   \put(0.2,1.8){$c$}
   \put(4.0,1.8){$\Longrightarrow$}       
   \put(1.0,0){$(\Sigma,\phi)$}
   \put(6.4,0){$(\Sigma',\phi')$}
  \end{picture}
  \caption{Stallings twist}
  \label{fig:twist}
 \end{center} }
 \end{figure}

 If an open book $(\Sigma, \phi)$ has a twisting loop on $\Sigma$,
 $(\pm 1)$-Dehn surgery along $c$ yield
 a new open book $(\Sigma', \phi')$ of $M$.
 We call this operation a {\itshape Stallings twist\/}
 along a twisting loop $c$.
 Note that the surface $\Sigma'$ is homeomorphic to $\Sigma$
 and its embedding into $M$ is changed (See Figure \ref{fig:twist}),
 and the monodromy map $\phi'= D(c)^{\pm 1} \circ \phi$,
 where maps act $\Sigma'$ from the left.

  A positive (negative resp.) stabilization of
  an open book $(\Sigma, \phi)$ of a closed oriented $3$-manifold
  is an open book $(\Sigma', \phi')$ of $M$ such that
  $\Sigma'$ is a plumbing of
  a positive Hopf band $H^{+}$ (negative Hopf band $H^{-}$ resp.)\
  and $\Sigma$.
  The new monodromy map $\phi'= D(\gamma)\circ\phi$
  ($\phi'= D(\gamma)^{-1}\circ\phi$ resp.),
  where $\gamma$ is the core curve of the Hopf band.
  we say that a stabilization along an arc properly embedded in $\Sigma$
  as a stabilization along a rectangle which is a regular
  neighbourhood of the arc in $\Sigma$.
 \begin{figure}[ht]
  {\unitlength=1cm
  \begin{center}
   \begin{picture}(12,3.5)(0,0)
    \put(0,0.6){\includegraphics[height=3cm]{%
    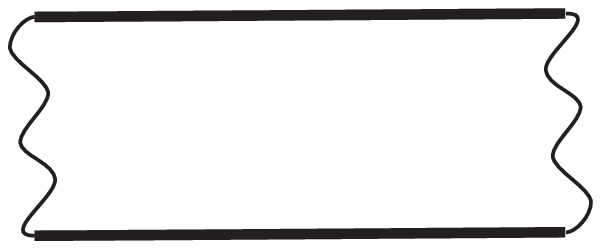} }       
    \put(4.0,0.6){\includegraphics[height=3cm]{%
    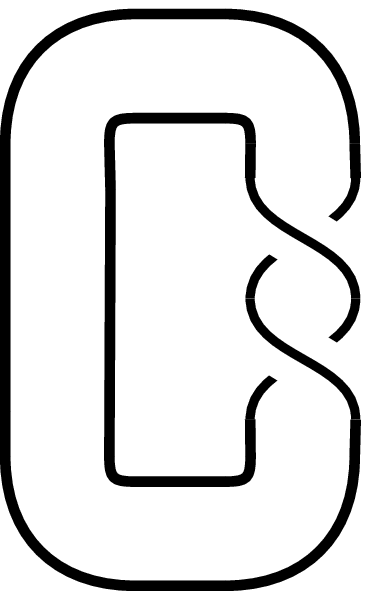} }
    \put(9,0.6){\includegraphics[height=3cm]{%
    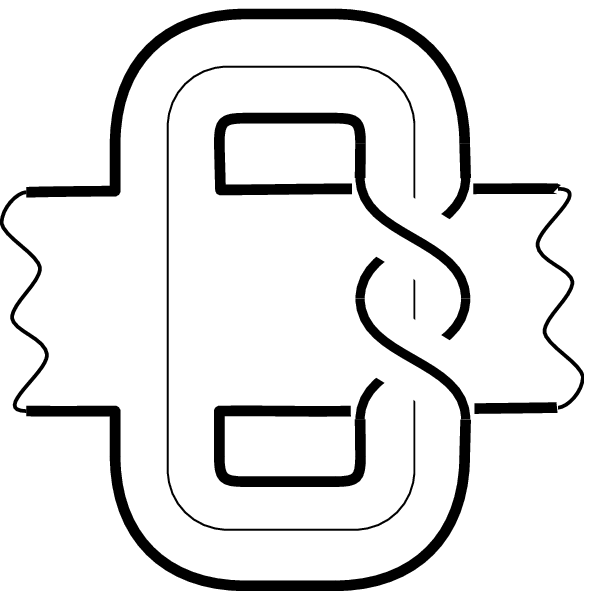} }
   \put(3.3,2.3){$\bigcup$}       
   \put(7.2,2.0){$\Longrightarrow$}       
   \put(1,0){$(\Sigma,\varphi)$}
   \put(5.8,.4){positive}
   \put(5.0,-0.2){Hopf band $H^{+}$}
   \put(9.6,2.0){$c$}
   \put(10.8,0){$(\Sigma',\varphi')$}
   \end{picture}
   \caption{positive stabilization of open book}
   \label{fig:stb^{+}}
  \end{center} }
 \end{figure}
 \subsection{Contact structure}
%
 A {\itshape contact form\/} on $M$ is
 a smooth global non-vanishing $1$-form $\alpha$ satisfying
 $\alpha\wedge d\alpha \neq 0$ at each point of $M$.
 A {\itshape contact structure\/} $\xi$ on $M$ is a $2$-plane field
 defined by the kernel of $\alpha$.
 The pair $(M,\xi)$ is called a {\itshape contact $3$-manifold\/}.
 We say that a contact structure
 $\xi=\operatorname{ker}\alpha$
 is {\itshape positive\/} when $\alpha\wedge d\alpha > 0$.
 We assume that a contact structure is positive throughout this article.

 We say that two contact structures on $M$, $\xi_{0}$ and $\xi_{1}$,
 are isotopic
 if there is a diffeomorphism $f: M \rightarrow M$ such that
 $\xi_{1}=f_{*}(\xi_{0})$.


 A simple closed curve $\gamma$ in contact $3$-manifold $(M,\xi)$
 is {\itshape Legendrian\/}
 if $\gamma$ is always tangent to $\xi$,
 i.e., for any point $x \in \gamma$,
 $T_x\gamma \subset \xi_x$.
 A Legendrian curve $\gamma$ has a natural framing
 called the {\itshape Legendrian framing\/}
 denoted by $\operatorname{Fr}(\gamma;\xi)$, 
 which is determined by a vector field on $\xi|_{\gamma}$
 such that each vector is transverse to $\gamma$.

 Let $E$ be a disk embedded in a contact manifold $(M,\xi)$.
 $E$ is an {\itshape overtwisted disk\/}
 if $\partial E$ is a Legendrian curve in $(M,\xi)$ and
 $\operatorname{Fr}(\partial E;E)
 =\operatorname{Fr}(\partial E;\xi)$.
%
%
%
 A contact structure $\xi$ on $M$ is {\itshape overtwisted\/}
 if there is an overtwisted disk $E$ in $(M,\xi)$.
 A contact structure is called {\itshape tight\/}
 if it is not overtwisted.

%
%

\subsection{Contact structures and open books}
%
 A contact structure $\xi$ on $M$ is said to be {\itshape supported\/}
 by an open book $(\Sigma,\phi)$
 if it is defined by a contact form $\alpha$
 such that
 (1) on each fiber $\Sigma_{t}$,
 $\ud{\alpha}|_{\Sigma_{t}} > 0$ and  
 (2) On $K=\partial\Sigma$,
 $\alpha(v_{p})>0$ for any point $p\in K$,
 where $v_{p}$ is a positive tangent vector of $K$ at $p$.
%
 W. Thurston and H. Winkelnkemper \cite{TW} showed that
 one can always construct a contact structure on $M$
 starting from a structure of an open book of $M$.
 The resulting contact structure is supported by the open book.

\begin{thm}[Giroux 2000,\cite{Gi};Torisu 2000,\cite{To}]\label{thm:2cs}
 Contact structures supported by the same open book are isotopic.   
\end{thm}

\begin{rmk}\label{rmk:close}
 It is known (e.g.\ \cite{Go}) that for an open book $(\Sigma,\phi)$
 there is a contact structure supported by $(\Sigma,\phi)$ such that
 at any point $p\in\operatorname{Int}\Sigma$
 the plain of $\xi$ is arbitrary close to the tangent plain
 of $\Sigma$.
 By Theorem \ref{thm:2cs},
 we may assume that a contact structure supported by an open book
 always has the property.
\end{rmk}

 As mentioned in Section \ref{sec:intro},
 Giroux showed that there is a one-to-one correspondence
 between contact structures and open books. 

\begin{thm}[Giroux 2000,\cite{Gi}] \label{thm:Giroux}
 Every contact structure of a closed oriented $3$-manifold is supported
 by some open books.
 Moreover open books supporting the same contact structure are
 equivalent up to positive stabilization. 
\end{thm}

 In Section \ref{sec:otob}
 we will need to show that
 a given simple closed curve on the fiber surface $\Sigma$
 of an open book $(\Sigma,\phi)$
 can be realized as a Legendrian curve
 in the contact structure supported by the open book.
 
 \begin{dfn}
  A simple closed curve $c$ on $\Sigma$ is {\itshape isolated\/}
  if there is a connected component $R$ of $\Sigma-c$ such that
  $R\cap \partial\Sigma=\emptyset$.
  We say that $c$ is {\itshape non-isolated\/} if it is not isolated.
 \end{dfn}
 
 The following lemma is a variant of the Legendrian Realization Principle
 on the convex surface theory,
 due to Ko Honda \cite{Ho}.
 \begin{lem}\label{lem:LeRP}
  There is a contact structure $\xi$ supported by $(\Sigma,\phi)$
  such that a simple closed curve $c$ on $\Sigma_{0}$ is Legendrian
  in $(M,\xi)$
  if and only if
  $c$ is non-isolated on $\Sigma_{0}$.
 \end{lem}
 \begin{proof}
  Let $(\Sigma,\phi)$ be an open book of $M$,
  and $\alpha$ a contact form on $M$
  which define a contact structure $\xi_{(\Sigma,\phi)}$
  supported by $(\Sigma,\phi)$.

  `only if' part.
  Let $c$ be a simple closed Legendrian curve in $\Sigma_{0}\subset(M,\xi)$.
  Suppose that $c$ is isolated in $\Sigma_{0}$,
  i.e., there is a subsurface $S\subset \Sigma_{0}$
  with $\partial S=c$.
  By the fact that $c$ is Legendrian and Stokes' theorem,
  we have that
  $0=\int_{c}\alpha = \int_{S}\ud{\alpha}$.
   This contradicts the fact that $\xi$ is supported by $(\Sigma,\phi)$,
  i.e., $\ud{\alpha}|_{\Sigma_{t}}>0$.
  Thus we have done a proof of `only if' part.

  `if' part.
  Let $c$ be a non-isolated simple closed curve on a fiber $\Sigma_{0}$
  of an open book $(\Sigma,\phi)$.
  We will construct a contact structure $\xi$ on $M$ supported by
  $(\Sigma,\phi)$, setting $c$ to be Legendrian.

  J.\ M.\ Montesinos and H.\ R.\ Morton \cite{MM} showed that
  for any open book $(\Sigma,\phi)$
  there are an open book $(D^{2},\beta)$ of $S^{3}$
  and simple covering map $\pi:M\rightarrow S^{3}$
  such that $\Sigma_{t}=\pi^{-1}(D_{t})$ for each $t\in[0,1)$
  and $\pi\circ\phi=\beta\circ\pi$,
  where $\{D_{t}\}$ is a family of fiber surfaces of $(D^{2},\beta)$.
  
  Since $c$ is non-isolated,
  $\Sigma$ has a handle decomposition such that
  $c$ is decomposed into two core arcs of $1$-handles and two arcs
  properly embedded in $0$-handles.
  Then it is easy to see that we can choose the covering map $\pi$
  so that $\pi(c)$ covers an arc $a$ on $D_{0}$.

  Take a contact form $\alpha$ on $S^{3}$
  which is supported by $(D^{2},\beta)$.
  We may assume that $\alpha|_{D_{0}}=x\ud{y}-y\ud{x}$
  in coordinates $(x,y)$ of $D$ such that
  the arc $a$ is contained in $x$-axis.
  Note that $x$-axis is Legendrian in
  $(S^{3},\operatorname{ker}(\alpha))$
  and so is $a$.
  Then the pullback $\pi^{*}\alpha$ is a contact form on $M$ supported
  by $(\Sigma,\phi)$ and $c$ is Legendrian
  in $(M,\operatorname{ker}(\pi^{*}\alpha))$.
 \end{proof}


%% file: otob.tex

 For the simplicity, we call an open book supporting an overtwisted
 contact structure an {\itshape overtwisted open book\/}.
 In this section,
 we first present two propositions which give
 sufficient conditions of an open books to be overtwisted,
 and then we give a proof of Theorem \ref{thm:main}
 using the propositions.

  \subsection{sufficient conditions for overtwistedness}


 Let $(\Sigma,\phi)$ be an open book of a closed oriented $3$-manifold
 $M$. 
   
  \begin{prop}
   \label{prop:twisting loop}
   $(\Sigma, \phi)$ is overtwisted
   if $(\Sigma, \phi)$ has a non-isolated twisting loop.
  \end{prop}
  \begin{proof}
   Let $c$ be a non-isolating twisting loop on $\Sigma$,
   and let $\xi_{(\Sigma, \phi)}$ denote a contact structure supported
   by $(\Sigma, \phi)$.
   By the definition of twisting loop,
   $c$ bounds a disk $D$ in $M$ such that
   \begin{equation}\label{eqn:01}
    \operatorname{Fr}(c;D)
     =\operatorname{Fr}(c;\Sigma).
   \end{equation}  
   Since $c$ is non-isolated in $\Sigma_{0}$,
   by Lemma \ref{lem:LeRP}
   we may assume that
   $c$ is Legendrian in $(M,\xi_{(\Sigma, \phi)})$.
   
   On the interior of $\Sigma$,
   plains of $\xi_{(\Sigma, \phi)}$ are arbitrary close to tangent
   plains of $\Sigma$ as mentioned Remark \ref{rmk:close}.
   So we have that
   \begin{equation}\label{eqn:02}
    \operatorname{Fr}(c;\xi_{(\Sigma, \phi)})
     =\operatorname{Fr}(c;\Sigma).
   \end{equation}
   From the equations (\ref{eqn:01}) and (\ref{eqn:02}) we have that
   \begin{equation*}\label{eqn:03}
    \operatorname{Fr}(c;\xi_{(\Sigma, \phi)})
     =\operatorname{Fr}(c;D).
   \end{equation*}
   This means that $D$ is an overtwisted disk in $(M,\xi_{(\Sigma,\phi)})$.
  \end{proof}


  Next we focus on an arc properly embedded on the fiber surface
  of an open book and its image of the monodromy map, 
  and show another criterion of overtwistedness of open books.
  
\begin{figure}[ht]
   {\unitlength=1cm
   \begin{center}
    \begin{picture}(10,4.6)(0,0)
     \put(0,.6){\includegraphics[height=3.6cm,keepaspectratio]{%
     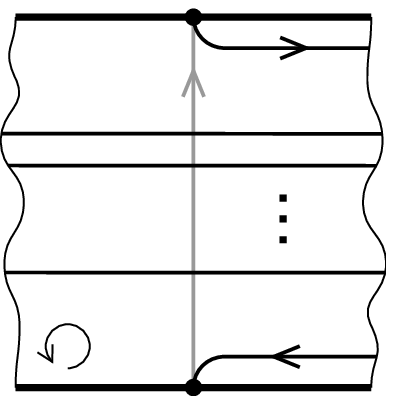}}
     \put(5.2,.7){\includegraphics[height=3.4cm,keepaspectratio]{%
     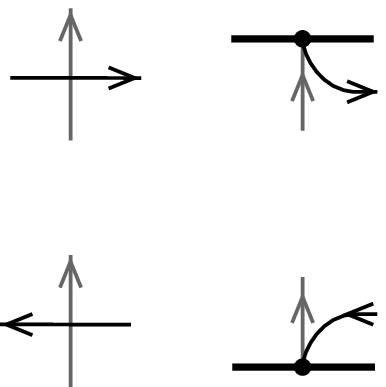}}
     \put(1.4,4.3){\textcircled{{\small $-$}}}
     \put(1.4,.2){\textcircled{{\small $+$}}}
     \put(1.5,1.0){$a$}
     \put(3.5,3.7){$\phi(a)$}
     \put(5.5,4.0){$a$}
     \put(6.2,3.1){$\phi(a)$}
     \put(5.5,1.8){$a$}
     \put(6.2,0.9){$\phi(a)$}
     \put(7.6,2.8){$a$}
     \put(8.5,3.1){$\phi(a)$}
     \put(7.6,1.7){$a$}
     \put(8.5,1.5){$\phi(a)$}
     \put(9.5,3.4){: $i_{p}=-1$}
     \put(9.5,1.3){: $i_{p}=+1$}
     \put(1.9,-.2){(a)}
     \put(7.5,-.2){(b)}
    \end{picture}
    \caption{}
    \label{fig:arc01}
   \end{center} }
  \end{figure}

  Let $a$ be an arc properly embedded in $\Sigma$.
  We always assume that $\phi(a)$ is isotoped relative to the boundary
  so that the number of intersection points between $a$ and $\phi(a)$
  is minimised.
 We orient the closed curve $a\cup\phi(a)$.
 It does not matter which orientation is chosen.
 At a point $p$ of $a\cap\phi(a)$ define $i_{p}$ to be $+1$
 if the oriented tangent to $a$ at $p$
 followed by the oriented tangent to $\phi(a)$ at $p$ is
 an oriented basis for $\Sigma$
 otherwise we set $i_{p}=-1$ (See Figure \ref{fig:arc01}).
 We define two kinds of intersection numbers of $a$ and $\phi(a)$
 as in Goodman's way \cite{Go};
 The geometric intersection number,
	   $i_{\mathrm{geom}}(a,\phi(a))
	   =\sum_{a\cap\phi(a)\cap\operatorname{Int}\Sigma}
	   |i_{p}|$,
	   is the number of intersection point of $a$ and $\phi(a)$
	   in the interior of $\Sigma$.
 The boundary intersection number,
	   $i_{\partial}(a,\phi(a))
	   =\frac{1}{2}
	   \sum_{a\cap\phi(a)\cap\partial\Sigma}
	   i_{p}$,
	   is one-half
	   the oriented sum over intersections at the boundaries of the arcs.

   \begin{prop}\label{prop:typeII}
    $(\Sigma, \phi)$ is overtwisted
    if $(\Sigma, \phi)$ has a proper arc $a$ such that 
    $a$ is not isotopic to $\phi(a)$ and satisfies
    \[
    i_{\mathrm{geom}}(a,\phi(a))
    =i_{\partial}(a,\phi(a))
    =0.
    \]
   \end{prop}
   \begin{proof}
    By the open book structure $(\Sigma,\phi)$ of $M$,
    we have a homeomorphism
    $h: E(\partial\Sigma)\rightarrow
    \Sigma\times [0,1]/(x,0)\sim(\phi(x),1)$
    with $h(\Sigma\cap E(\partial\Sigma))=\Sigma_{0}=\Sigma_{1}$.
    Put $\Delta=h^{-1}(h(a)\times [0,1])$.
    Let $N$ be a small neighbourhood of $a$ on $\Sigma$,
    and $a'$ a parallel copy of $a$ in $N$
    such that $a'\cap \phi(a)$ is two points, say $p$ and $q$.

    Deform $\Delta$ as follows;
    Push small neighbourhoods of $p$ and $q$ on $\Delta$ along $a'$
    to make them close,
    and merge them so that we have a half pipe along $a'$.
    Now we obtain an annulus $\Delta'$
    such that one of the boundary components is a
    simple closed curve on $\Sigma$ isotopic to $a\cup\phi(a)$
    and another is consists of $a$ and an arc $a''$
    isotopic to $a$ relative to $\partial a$.
    Slide $\Delta'$ on $\Sigma$ so that $a''$ overlap with $a$.
    Finally we see that
    $\Delta'$ becomes a disk $D$ embedded in $M$ such that
    $\partial D$ is a simple close curve on $\Sigma$
    and $D\cap \Sigma=a$(See Figure \ref{fig:deform}).
    
    It follows from the construction that $\partial D$ is a twisting loop
    on $\Sigma$.
    \begin{figure}[ht]
     {\unitlength=1cm
     \begin{center}
      \begin{picture}(10,2.6)(0,0)
       \put(0,-.5){\includegraphics[height=3.2cm]{%
       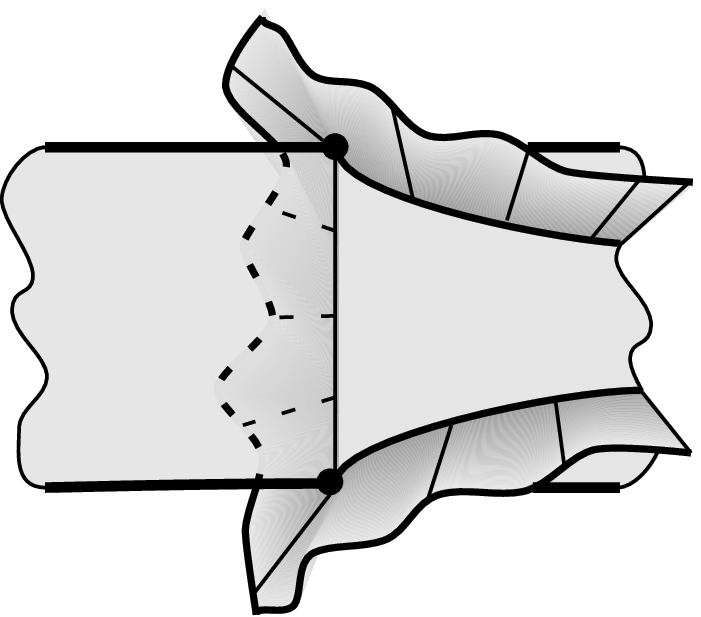} }       
       \put(6,-.5){\includegraphics[height=3.2cm]{%
       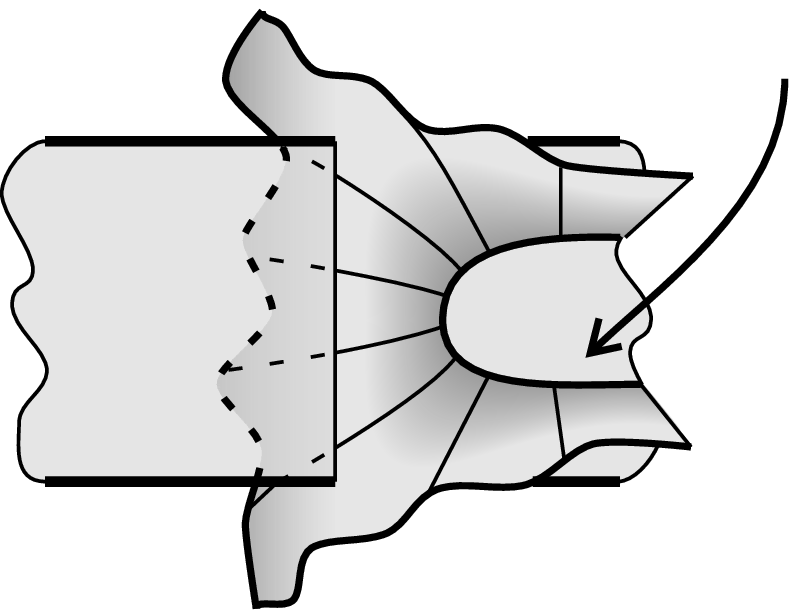} }
       \put(1.9,1.1){$a$}
       \put(3.5,1.3){$\phi(a)$}
       \put(4.6,1.1){$\xrightarrow[\text{deform}]{}$}
       \put(9,2.7){a twisting loop}
      \end{picture}
      \caption{}
      \label{fig:deform}
     \end{center} }
    \end{figure}

    If $\partial D$ is non-isolated in $\Sigma$,
    by Proposition \ref{prop:twisting loop}
    we have that $D$ is an overtwisted disk in
    $(M,\xi_{(\Sigma,\phi)})$,
    where $\xi_{(\Sigma,\phi)}$ is a contact structure supported by
    $(\Sigma,\phi)$.

    Suppose that $\partial D$ is isolated.
    Note that $\partial D$ is isotopic to $a\cup\phi(a)$.
    So we have a connected component $S$ of $\Sigma-(a\cup \phi(a))$
    such that $S\cap \partial\Sigma=\emptyset$.
    We will show that 
    we can obtain a new open book $(\Sigma',\phi')$
    by positive stabilizations such that
    $a$ and $\phi'(a)$ satisfy the assumption of this proposition
    and $a\cup\phi'(a)$ is non-isolated in $\Sigma'$. 

    Note that $S$ has the genus greater than $0$, since $a$ and
    $\phi(a)$ is not isotopic relative to the boundary.
    Let $\beta_{1}$ be an arc properly embedded in $\Sigma$
    as shown in Figure \ref{fig:ToBeNon-isolated} (right)
    such that 
    the arc $\beta_{1}\cup S$ is not boundary-parallel
    and non-separating in $S$.
    By a positive stabilization along $\beta_{1}$
    we have an open book $(\Sigma'',\phi'')$ and 
    we can find an arc $\beta_{2}$ in $\Sigma''$ such that
    $\beta_{2}$ intersects with $a$ and $\phi''(a)$ at just one point
    each as shown in Figure \ref{fig:ToBeNon-isolated} (center). 
    Then a positive stabilization along $\beta_{2}$ yields
    an open book $(\Sigma',\phi')$,
    and we can easily see that $a\cup\phi'(a)$ is non-isolated on
    $\Sigma'$ (Figure \ref{fig:ToBeNon-isolated}).
\begin{figure}[h]
  {\unitlength=1mm
  \begin{picture}(120,26)(0,0)
   \put(0,0){\includegraphics[width=2.8cm,keepaspectratio]{%
   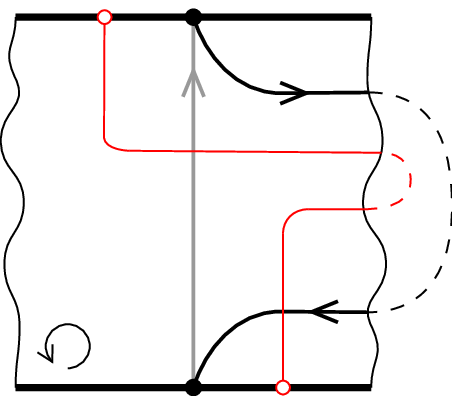}}
   \put(46,-4){\includegraphics[width=2.8cm,keepaspectratio]{%
   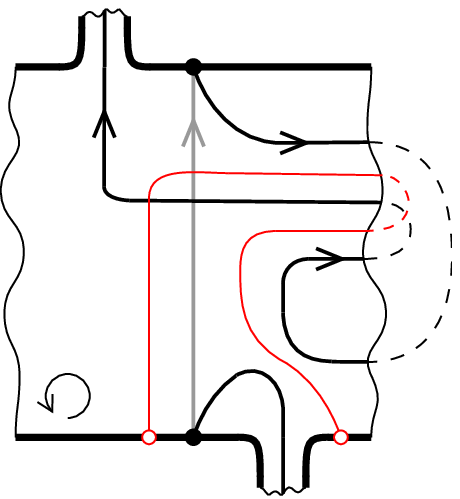}}
   \put(92,-4){\includegraphics[width=2.8cm,keepaspectratio]{%
   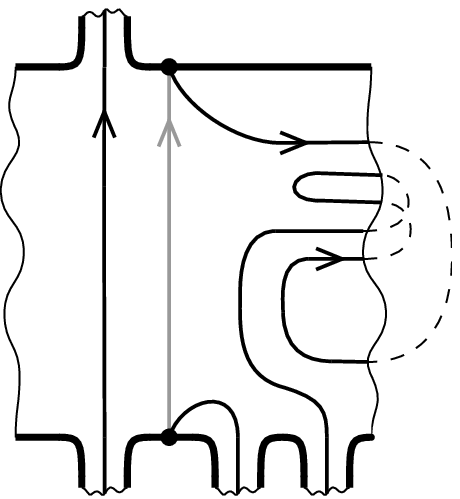}}
   \put(30,10){$\xrightarrow[\text{along } \beta_{1}]{\operatorname{stb}^{+}}$}
   \put(76,10){$\xrightarrow[\text{along } \beta_{2}]{\operatorname{stb}^{+}}$}
   \put(2,17){$\small\beta_{1}$}
   \put(51,7){$\small\beta_{2}$}
%
  \end{picture} }
 \caption{}
 \label{fig:ToBeNon-isolated}
\end{figure}     
   \end{proof}

  \subsection{Proof of Theorem \ref{thm:main}}

  \noindent\underline{(1) $\Rightarrow$ (2)}
  Let $D$ be an overtwisted disk embedded in $(M,\xi_{(\Sigma,\phi)})$.
  By Corollary 4.23 in \cite{Et2},
  we can construct an open book $(\Sigma_{0},\phi_{0})$
  supporting $\xi_{(\Sigma,\phi)}$
  such that $\partial D$ is on $\Sigma_{0}$.
  It is easy to see that $\partial D$ is a twisting loop on
  $\Sigma_{0}$.
  Since both $(\Sigma,\phi)$ and $(\Sigma_{0},\phi_{0})$ support the
  same contact structure,
  Giroux's one-to-one correspondence tells us that
  they are equivalent up to positive stabilization.
  
  \medskip
  \noindent\underline{(2) $\Rightarrow$ (1)}
  This part immediately follows from Proposition \ref{prop:twisting loop}
  and Giroux's theorem.

  
  \medskip
  \noindent\underline{(1) $\Rightarrow$ (3)}
  Goodman showed in \cite[Theorem 5.1]{Go} that an overtwisted contact
  structure on a closed oriented $3$-manifold has a supporting open book
  with a sobering arc,
  and a sobering arc has the boundary points satisfing $i_{\partial}\leq 0$.
  
\begin{figure}[ht]
    {\unitlength=1cm
    \begin{center}
     \begin{picture}(10,3.0)(0,0)
      \put(0,.3){\includegraphics[width=2.4cm]{%
      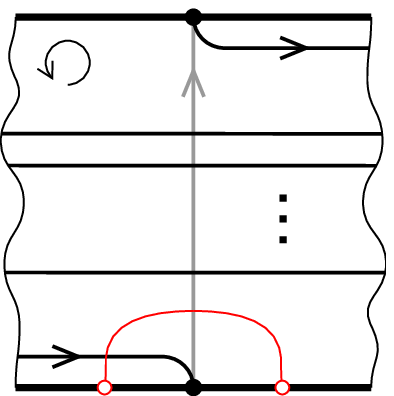} }
      \put(6.5,-.2){\includegraphics[width=2.4cm]{%
      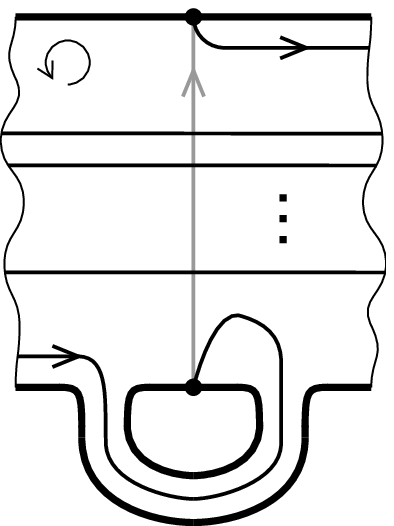} }
      %
      %
      \put(3.5,1.3){
      $\xrightarrow{\operatorname{stb}^{+} \text{along } \beta}$ }
      \put(.9,1.3){$a$}
      \put(1.8,0.6){$\beta$}
      %
      \put(.8,2.8){\textcircled{{\small $-$}}}
      \put(.8,-.0){\textcircled{{\small $-$}}}
      %
      \put(1.2,2.8){$x_{0}$}
      \put(1.2,-.0){$x_{g+1}$}
      \put(2.5,0.2){$(\Sigma,\phi)$}
      \put(8.8,0.2){$(\Sigma',\phi')$}
      %
     \end{picture}
    \caption{We may assume that $i_{\partial}=0$.}
    \label{fig:stab-00}
    \end{center} }\end{figure}

  \medskip
  \noindent\underline{(3) $\Rightarrow$ (1)}
  Let $a$ be a properly embedded arc on $\Sigma$ with
  $i_{\partial}(a,\phi(a))\leq 0$,
  The arc $a$ has at least one negative endpoint, say $x_{0}$.
  We orient $a$ so that a oriented tangent vector
  of $a$ at $x_{0}$ is outward from $\Sigma$.
  Put $g = i_{\mathrm geom}(a,\phi(a))$.
%
  %
  Stating from $x_{0}$,
  we assign $x_{1},\dots,x_{g+1}$ to the points of $a \cap \phi(a)$.
  Suppose that the another endpoint $x_{g+1}$ is also negative,
  i.e., $i_{\partial}(a,\phi(a))=-1$.
  Let $\beta$ be a small properly embedded arc on $\Sigma$ rounding
  $x_{g+1}$. 
  By a positive stabilization along $\beta$
  we have a new open book with a new monodromy map,
  on which the point $x_{g+1}$ is positive
  (See Figure \ref{fig:stab-00}).
  Thus we may assume that $i_{\partial}(a,\phi(a))=0$.
    %
  

  We show this part by induction on $g$.
  In the case where $g=0$,
  we have already proven as Proposition \ref{prop:typeII}
  that we can find an overtwisted disk in $(M,\xi_{(\Sigma,\phi)})$
  with the boundary on $\Sigma_{0}$.

\begin{figure}[ht]
   {\unitlength=1cm
   \begin{center}
    \begin{picture}(5.6,4.5)(0,0)
     \put(0,.2){\includegraphics[height=4.0cm]{%
     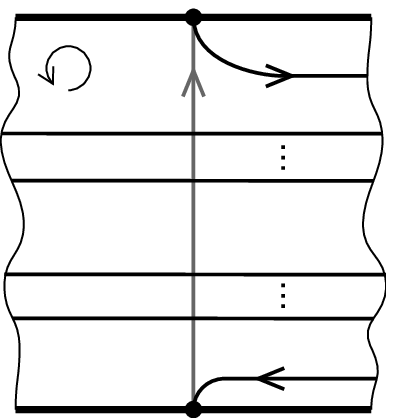} }
     %
     %
     \put(1.4,4.4){$x_{0}$}   
     \put(1.4,3.2){$x_{1}$}   
     \put(1.0,2.3){$x_{k-1}$}   
     \put(1.4,1.8){$x_{k}$}   
     \put(1.4,.9){$x_{g}$}   
     \put(1.0,-.1){$x_{g+1}$}   
     \put(2.8,3.7){$R_{0}$}   
     \put(2.8,3.15){$R_{1}$}   
     \put(2.8,1.8){$R_{k}$}   
     \put(2.8,.75){$R_{g+1}$}   
     \put(1.9,4.4){\textcircled{{\small $-$}}}
     \put(1.9,-.2){\textcircled{{\small $+$}}}
     \put(2.0,2.0){$a$}
     \put(3.7,3.5){$\phi(a)$}
     \put(4.0,0){$\operatorname{nbd}(a;\Sigma)$}
    \end{picture}
    \caption{}
    \label{fig:naming}
   \end{center} }\end{figure}

  Suppose that $g>0$.
%
  Let $\alpha_{i}$ be a sub-arc of $a$ connecting $x_{i-1}$ and $x_{i}$
  for $1\leq i\leq g+1$ and 
  $\gamma_{i}$ a connected component of $N \cap \phi(a)$
  containing the point $x_{i}$ for $0\leq i\leq g+1$.
  We denote by $N_{R}$ a connected component of $N-a$ which has
  intersection with $\gamma_{0}$.
  Let $R_{0}$ be a connected region of $N_{R}-\gamma_{0}$ such that
  $R_{0}\cap a=x_{0}$,
  and $R_{i}$ a connected region of $N_{R}-\phi(a)$
  such that $R_{i}\cap a=\alpha_{i}$ for $1\leq i\leq g+1$ 
  (See Figure \ref{fig:naming}).
  We denote by $\widehat{R_{k}}$ a connected region of
  $\Sigma-(a\cup\phi(a))$ containing $R_{k}$.
  Note that some regions of
  $\widehat{R_{0}},\widehat{R_{1}},\dots,\widehat{R_{g+1}}$
  might indicates the same one.

 Tracing $\phi(a)$ along its orientation
 and picking up $x_{i}$'s on the points of $a\cap\phi(a)$,
 we obtain a word
 $w=x_{0}^{-1}x_{p(1)}^{\epsilon_{p(1)}}x_{p(2)}^{\epsilon_{p(2)}}
 \dots x_{p(g)}^{\epsilon_{p(g)}}x_{g+1}$,
 where $p$ is a permutation of $\{1,2,\dots,g\}$
 and $\epsilon_{k}$ is the sign of the point $x_{k}$ for $1\leq k \leq g$.
   %



 Now we introduce two types of positive stabilization which 
 reduces $g$, keeping $i_{\partial}=0$.
 
\begin{figure}[ht]
    {\unitlength=1cm
    \begin{center}
     \begin{picture}(12,3.3)(0,0)
      \put(0,.3){\includegraphics[width=4cm]{%
      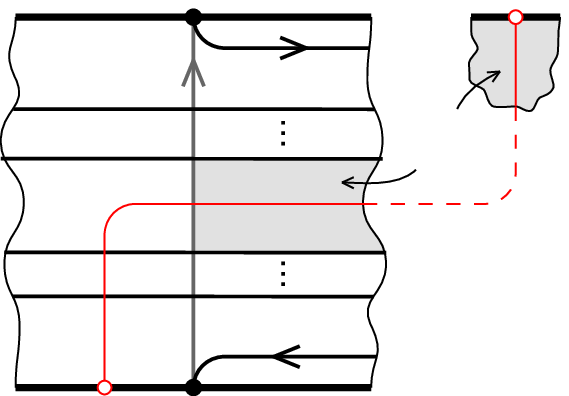} }
      \put(8.0,-.1){\includegraphics[width=4cm]{%
      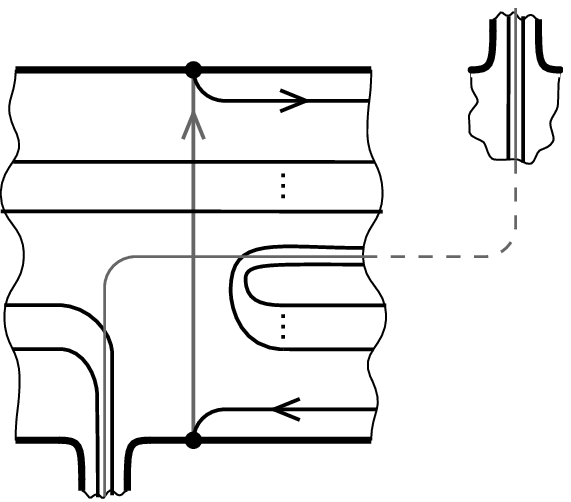} }
      %
      %
      \put(5.0,1.4){$\xrightarrow{\operatorname{stb}^{+}
      \text{along } \beta}$ }
      \put(3.0,1.95){$\widehat{R_{k}}$}
      \put(1.1,.6){$a$}
      \put(9.1,.6){$a$}
      \put(3.8,1.6){$\beta$}
      \put(11.8,1.6){$\hat{\beta}$}
      \put(1.2,3.2){\textcircled{{\small $-$}}}
      \put(1.2,0){\textcircled{{\small $+$}}}
      \put(9.2,3.2){\textcircled{{\small $-$}}}
      \put(9.4,0){\textcircled{{\small $+$}}}
      \put(3.2,0.2){$(\Sigma,\phi)$}
      \put(9.8,-0.4){$(\Sigma',\phi'=D(\hat{\beta}) \circ \phi)$}
     \end{picture}
     \caption{Reduction A}
     \label{fig:Red-A}
    \end{center} }\end{figure}

 First one is as follows.
 Suppose that there exist an integer $k$ ($1 \leq k \leq g$)
   such that 
   $\widehat{R_{k}}\cap \partial \Sigma\neq\emptyset$.
   Then there is an arc $\beta$ properly embedded in $\Sigma$
   as shown in Figure \ref{fig:Red-A} (left).
   We obtain an open book $(\Sigma',\phi')$
   by positive stabilization along $\beta$
   such that
   \[
   i_{\mathrm geom}(a,\phi'(a))= k-1 < g.
   \]
   We call this type of stabilization
   {\itshape reduction A\/}.

\begin{figure}[ht]
 {\unitlength=1cm
  \begin{center}
  \begin{picture}(11.5,3.3)(0,0)
   \put(0,0){\includegraphics[width=3.6cm, keepaspectratio]{%
   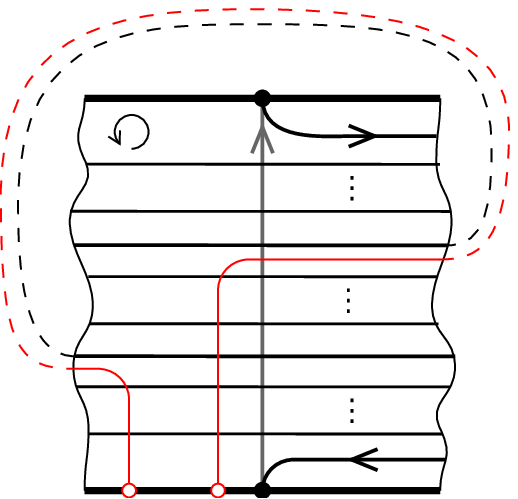}}
   \put(8,-.6){\includegraphics[width=3.6cm, keepaspectratio]{%
   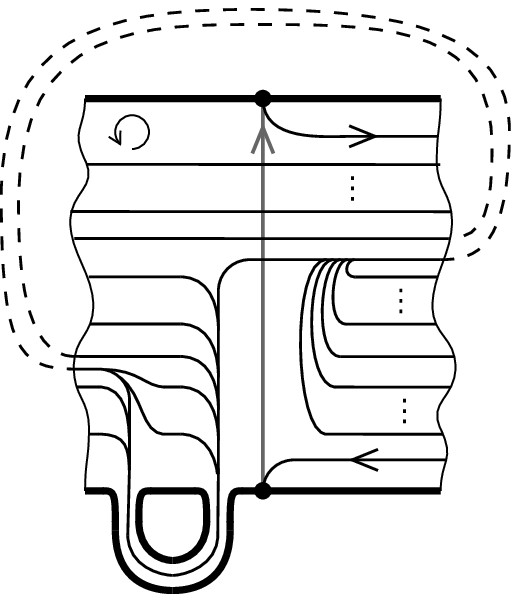} }
   \put(3.6,1.8){$\beta$}
   \put(5.0,1.4){$\xrightarrow{\operatorname{stb}^{+}
   \text{along } \beta}$ }
   %
%
%
%
  \end{picture}
   \caption{Reduction B}
   \label{fig:Red-B}
  \end{center}}\end{figure}

 Next,
 suppose that there are points $x_{i}$ and $x_{j}$ ($1\leq i<j \leq g$)
 in $a\cap\phi(a)$ such that
 they are adjacent in $\phi(a)$,
 their signs are the same,
 and the order $x_{i}, x_{j}$ is agree with
 the orientation of $\phi(a)$ if they have positive sign,
 not agree if negative.
 Then we can find an arc $\beta$ properly embedded in $\Sigma$
 as shown in Figure \ref{fig:Red-B} (left).
 We obtain an open book $(\Sigma',\phi')$
 by positive stabilization along $\beta$
 such that
 \[
 i_{\mathrm geom}(a,\phi'(a))= g-(j-i) < g.
 \]
 We call this type of stabilization
 {\itshape reduction B}.
%


 We may assume that
 there are no proper arcs which admit reduction A or reduction B.
 Then the neighbourhood of $a$ in $\Sigma$ has the following properties;
 \begin{enumerate}
  \item $\widehat{R_{k}}\cap \partial \Sigma=\emptyset$
	for all $k=1,2,\dots,g$,
  \item the word $w$ contains no consecutive letters 
	as $x_{i}^{-1}x_{j}^{-1}$
	or $x_{j}x_{i}$,
	where $0\leq i < j \leq g$.
 \end{enumerate}

 Now we focus on the point $x_{p(1)}$.
 Suppose that $\epsilon_{p(1)}=-1$,
 i.e., the first two letters of $w$ is $x_{0}^{-1}x_{p(1)}^{-1}$.
 This contradicts the property $(2)$ above.
 Thus we have that $\epsilon_{p(1)}=+1$.
%
%
 Suppose that $p(1)<g$.
 We can easily see that
 $\widehat{R_{p(1)+1}}\cap\partial\Sigma\neq\emptyset$.
 This contradicts the property $(1)$.
 Thus $p(1)=g$.
\begin{figure}[h]
  {\unitlength=1mm
  \begin{picture}(30,25)(0,0)
   \put(0,0){\includegraphics[height=2.5cm,keepaspectratio]{%
   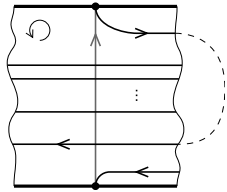}}
 %
 %
 %
  \end{picture} }
 \caption{$w=x_{0}^{-1}x_{g}\cdots x_{g+1}$}
 \label{fig:p(1)=g}
\end{figure}

\begin{figure}[h]
  {\unitlength=1mm
  \begin{picture}(100,33)(0,0)
   \put(0,0){\includegraphics[width=3.6cm,keepaspectratio]{%
   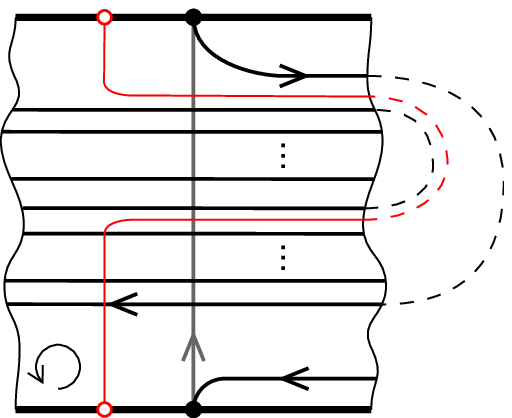}}
   \put(40,14){$\xrightarrow{\operatorname{stb}^{+}
   \text{along } \beta}$ }
   \put(65,-4){\includegraphics[width=3.6cm,keepaspectratio]{%
   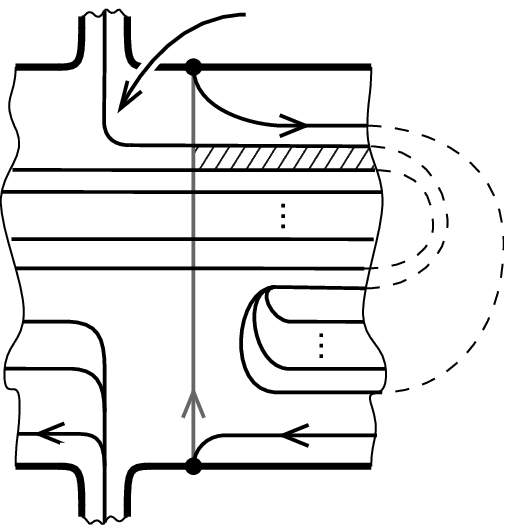} }
   \put(83,31){$g-k$ parallel arcs}
  \end{picture} }
 \caption{}
 \label{fig:Case1}
\end{figure}

 \noindent\underline{
 {\bfseries Case 1.}~~ $g > 1$.}
 We look at the letter $x_{1}$ in the word $w$.
 By the property $(2)$,
 we have that there is an integer $k$ ($2\leq k \leq g-1$) such that
 $x_{k}$ is adjacent to $x_{1}$ in $w$
 and they appear as
 $x_{1}^{-1}x_{k}$ or $x_{k}^{-1}x_{1}$.
 Put $R=\bigcup_{i=1}^{g}R_{i}$.
 We can find an arc $\beta$ properly embedded in $\Sigma$
 as shown in Figure \ref{fig:Case1} (left)
 such that $\beta\cap R$ is an arc not boundary-parallel in $R$
 and $\beta$ does not intersect with $\phi(a)$ in $R$.
 The positive stabilization along $\beta$,
 which keeps the intersection number of $a$
 with its image of the monodromy map, 
 yields a new open book, say $(\Sigma',\phi')$
 (See Figure \ref{fig:Case1}).
 Assign the names of regions $R_{0},R_{1},\dots,R_{g+1}$ to the new
 regions of $N-(a\cup\phi'(a))$ in the same manner.
 It is easy to see that the region $R_{g-k+1}$
 (shaded in Figure \ref{fig:Case1} (right))
 and $R_{0}$ are connected in $\Sigma'$, i.e.,
 $\widehat{R_{g-k+1}}=\widehat{R_{0}}$.
 Thus we can perform reduction A.
%

\begin{figure}[h]
  {\unitlength=1mm
  \begin{picture}(100,33)(0,0)
   \put(0,-2){\includegraphics[width=3.4cm,keepaspectratio]{%
   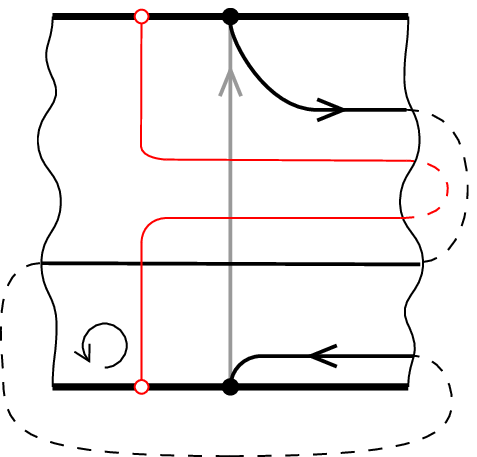}}
   \put(40,16){$\xrightarrow{\operatorname{stb}^{+}
   \text{along } \beta}$ }
   \put(65,-2){\includegraphics[width=3.4cm,keepaspectratio]{%
   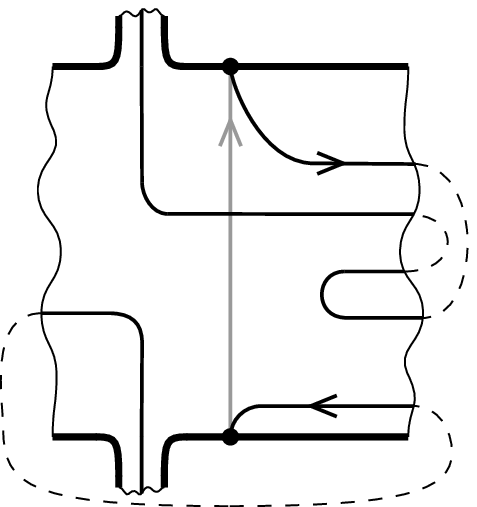} }
   \put(7,25){$\small\beta$}
%
  \end{picture} }
 \caption{}
 \label{fig:Case2}
\end{figure}

 \noindent\underline{
 {\bfseries Case 2.}~~ $g = 1$.}
 Put $R=\widehat{R_{1}}$.
 Recall that $R$ has no intersection with $\partial\Sigma$
 and the genus of $R$ is greater than $0$.
 We have an arc $\beta$ properly embedded in $\Sigma$
 as shown in Figure \ref{fig:Case2} such that
 $\beta\cap R$ is an arc not boundary-parallel, non-separating in $R$.
 
 As in the previous case,
 the positive stabilization along $\beta$ keeps
 the intersection number of $a$ with its image of the monodromy map.
 Let $(\Sigma',\phi')$ be a resulting open book and
 reassign the names of regions $R_{0},R_{1},R_{2}$ to the new
 regions of $N-(a\cup\phi'(a))$ in the same manner.
 Since $\beta$ is non-separating in $R$,
 the regions $R_{1}$ and $R_{2}$ are in the same connected
 component of $\Sigma-(a\cup\phi(a))$.
 Then $R_{1}$ and $R_{0}$ are also. 
 Thus we can perform reduction A.
%
 \qed
%
    